\documentclass{elsarticle}

\usepackage{graphicx}
\usepackage{latexsym}
\usepackage{amsmath,amsfonts}
\usepackage{color}
\usepackage{subfig}
\usepackage{multirow}

\newtheorem{theorem}{Theorem}[section]

\newtheorem{proposition}[theorem]{Proposition}
\newtheorem{corollary}[theorem]{Corollary}

\newenvironment{proof}{\noindent {\it Proof.}}{$\Box$}
\newenvironment{remark}{\noindent {\bf Remark.}}{}

\newcommand{\comp}[1]{\overline{#1}}
\newcommand{\neig}{N}
\newcommand{\ane}[1]{\bar \neig(#1)}
\newcommand{\pane}[1]{\bar \neig^{\textsc{u}}(#1)}
\newcommand{\nane}[1]{\bar \neig^{\textsc{l}}(#1)}
\newcommand{\REP}[1]{\REPORD{\prec}{#1}}
\newcommand{\REPORD}[2]{\mathcal{R}^{#1}_{#2}}
\newcommand{\COL}[1]{\COLORD{\prec}{#1}}
\newcommand{\COLORD}[2]{P^{\textsc{Co}}_{#1}(#2)}

\newcommand{\PREEXT}[1]{\PREEXTORD{\prec}{#1}}
\newcommand{\PREEXTORD}[2]{P^{\textsc{Pr}}_{#1}(#2)}
\newcommand{\PC}{\rho}

\newcommand{\X}[1]{\textsc{x}_{[#1]}}

\newcommand{\R}{\mathbb{R}}

\newcommand{\N}{\mathbb{N}}

\newcommand{\comment}[1]{{\color{red} {[#1]}}}
\newcommand{\edit}[1]{{\color{blue} {#1}}}

\newcommand{\kite}{kite}

\usepackage{array,booktabs}

\begin{document}

\begin{frontmatter}

\title{Polyhedral studies of vertex coloring problems:\\The asymmetric representatives formulation}


\author{Victor Campos\fnref{a}}
\author{Ricardo C. Corr\^ea\fnref{ufrrj,fromufc}}
\author{Diego Delle Donne\fnref{b,c,d}}
\author{\\Javier Marenco\fnref{b,c}}
\author{Annegret Wagler\fnref{d}}

\fntext[fromufc]{This author was with Universidade Federal do
	Cear\'a, Dto. de Computa\c c\~ao, Brazil, when part of this work was done.}

\address[a]{Federal University of Cear\'a, Fortaleza, Brazil}
\address[ufrrj]{Universidade Federal Rural do Rio de Janeiro,
	Dto. de Ci\^encia da Computa\c c\~ao, RJ, Brazil}
\address[b]{Computer Sc. Department, FCEN, University of Buenos Aires, Buenos Aires, Argentina}
\address[c]{Sciences Institute, National University of General Sarmiento, Buenos Aires, Argentina}
\address[d]{LIMOS, Blaise Pascal University, Clermont Ferrand, France}

\date{}

\begin{abstract}
Despite the fact that some vertex coloring problems are polynomially solvable on certain graph classes, most of these problems are not ``under control'' from a polyhedral point of view.
The equivalence between \emph{optimization} and \emph{polyhedral separation} suggests that, for these problems, there must exist formulations admitting some elegant characterization for the polytopes associated to them. 
Therefore, it is interesting to study known formulations for vertex coloring with the goal of finding such characterizations.
In this work we study the asymmetric representatives formulation and we show that the corresponding coloring polytope, for a given graph $G$, can be interpreted as the stable set polytope of another graph $\REPORD{} G$, obtained from $G$. 
This result allows us to derive complete characterizations for the corresponding coloring polytope for some families of graphs, based on known complete characterizations for the stable set polytope. 

\noindent\textsc{Keywords}: vertex coloring, representatives formulation, polyhedral characterization.
\end{abstract}

\end{frontmatter}

\section{Introduction}
	Given a graph $G=(V,E)$, a \emph{coloring} of $G$ is an assignment $c:V\to\N$ of ``colors'' to the vertices of $G$ such that $c(v)\ne c(w)$ for each edge $vw \in E$. The \emph{vertex coloring problem} consists in finding a coloring of $G$ minimizing the number of used colors. This parameter is the \emph{chromatic number} of $G$, and is denoted by $\chi(G)$. 
	There exist in the literature many variants of the classical vertex coloring problem. 
	For instance \emph{$\mu$-coloring} \cite{B-C-mp-graco} and \emph{$(\gamma,\mu)$-coloring} \cite{MARENCO06} are such generalizations. 
	These problems take functions \mbox{$\gamma,\mu:V\to\mathbb{N}$} defining lower and upper bounds of the color to be assigned to each vertex, i.e., the obtained coloring $c:V\to\mathbb{N}$ must satisfy $\gamma(v) \leq c(v) \le \mu(v)$, for every $v \in V$ (in $\mu$-coloring, $\gamma$ is omitted).
	Another interesting variant (which generalizes both variants mentioned above) is the \emph{list coloring problem} \cite{Tuza}, which considers a set $L(v)$ of valid colors for each $v\in V$ and asks for a coloring $c$ such that $c(v)\in L(v)$ for all $v\in V$.
	This problem is one of the most general versions of coloring problems.
	Constraints coming from real-life settings motivate many other variations of the classical graph coloring problem.
	In this work, we focus on the following two generalizations:
\begin{description}
	\item{\em Precoloring extension} \cite{B-H-T-mp}: Given a graph $G=(V,E)$ and a partial assignment \mbox{$\rho:V'\to\N$}, for some $V'\subseteq V$, find a coloring $c$ with the smallest number of used colors such that $c(v) = \rho(v)$ for every vertex $v\in V'$. In other words, a subset $V'$ of vertices from $G$ is already colored and the problem asks to extend this coloring to the whole set $V$ in a minimum fashion.
	
	\item {\em Max-coloring} \cite{Demange02,Pemmaraju04}: Given a graph $G=(V,E)$ and a weight function \mbox{$\omega: V\to\R_+$}, find a proper coloring $c:V\to\N$ of $G$ minimizing the sum of the weights over all color classes, where the weight of a color class $S$ is given by the maximum weight of a vertex from $S$.
\end{description}
	
Although the classical vertex coloring problem is NP-hard \cite{Garey}, there are several graph classes for which this problem can be solved in polynomial time, 
most notably
perfect graphs \cite{Golumbic}. 
A graph $G$ is said to be \emph{perfect} if $\chi(H) = \omega(H)$ for every induced subgraph $H$ of $G$, where $\omega(H)$ represents the size of a maximum clique of $H$. 
However, the variants of the coloring problem mentioned above may not be polynomially solvable for perfect graphs and for other graph classes where the classical vertex coloring problem is indeed.
For instance, in \cite{MARENCO06,BONOMOUNP}, the complexity boundary between coloring and list-coloring is studied for several subclasses of perfect graphs. 
It is shown that the precoloring extension problem can be solved in polynomial time for cographs, complete bipartite, split and complete split graphs, among others. On the other side there are many other classes for which this problem is NP-hard (see \cite{MARENCO06,BONOMOUNP} for details).
The max-coloring problem is also NP-hard in general
and
it seems to be substantially harder than the classical vertex coloring problem; for instance it is NP-hard in chordal graphs \cite{Demange02} and in interval graphs \cite{Pemmaraju04}, although in both graph classes, the classical version has linear resolution \cite{Golumbic}.\\

\emph{Integer linear programming} (ILP) has proved to be a very suitable tool for solving combinatorial optimization problems \cite{N-W-book}, and in the last decade ILP has been successfully applied to graph coloring problems, by resorting to several formulations for the classical vertex coloring problem. 
The {\em standard model} \cite{Coll, MZ06, MZ08} includes a binary variable $x_{ic}$ for each vertex $i\in V$ and each color $c\in C$, where $C$ represents the set of available colors, asserting whether vertex $i$ is assigned color $c$ or not. This formulation may be extended with variables $w_c$ for each color $c \in C$ specifying whether this color is used or not; the minimum coloring is found by minimizing the sum of the latter variables.\\

There are many other interesting formulations for the vertex coloring problem and it variants, namely: the orientation model \cite{BoEiGrMa98b}, the distance model \cite{DelleDonne09}, the supernodal formulation \cite{Burke} and the well-known column generation approach by Mehrotra and Trick \cite{Mehrotra}. 
In this work, we focus on the following formulation:

\begin{description}
	

	\item {\em Asymmetric representatives formulation} \cite{Campelo08,Campelo04}: In this formulation, a coloring is determined by the color classes it induces, and each class is represented by one of its members. For each non-adjacent pair of vertices $i,j\in V, i < j$, the model uses a binary variable $x_{ij}$ stating whether vertex $i$ is the ``representative'' of the color class assigned to vertex $j$ or not. Additionally, the model uses the variable $x_{ii}$ for each $i\in V$, asserting whether $i$ is the representative of its own color class or not. 
	We describe this model in detail in the following section.

\end{description}

Although ILP is NP-hard, in many cases a complete description of the convex hull of its solutions is known and this description can be used to solve the \emph{separation problem} associated to this polytope in polynomial time \cite{SCH03}. 
Based on the ellipsoid method, Gr\"otschel, Lov\'asz and Schrijver \cite{GLS88} proved that the separation problem and the optimization problem over a polytope are polynomially equivalent, i.e., if one problem is polynomially solvable, so is the other one (we properly define these two problems in the following section). From this equivalence stems the general belief that if a combinatorial optimization problem can be solved in polynomial time, then there should exist some ILP formulation of the problem for which the convex hull of its solutions admits an ``elegant'' characterization.\\

Despite the fact that some vertex coloring problems are polynomially solvable on certain graph classes, most of these problems are not ``under control'' from a polyhedral point of view. 
The mentioned equivalence between optimization and separation suggests that, for these
problems, there must exist formulations with polynomially solvable separation problems and,
moreover, that these formulations may admit elegant characterizations. 
The search for
such characterizations is the main objective and motivation of our work.\\

From a theoretical point of view, our main objective is to complete the polyhedral counterpart of these combinatorially solved graph coloring problems. On the other side, the study of these polytopes may lead us to a better understanding of their structures allowing us to (polyhedrally) find new classes of graphs colorable in polynomial time.\\

In this sense, some work has been done on the standard formulation in a previous work \cite{DDDMarencoDO} and nice characterizations were found for simple graph classes such as trees and block graphs. 
Also, the separation problem associated to this formulation for the classical vertex coloring problem cannot be solved in polynomial time for graph classes for which list-coloring is NP-hard, unless P=NP, meaning that this formulation falls short too soon. Therefore, one direct next step on this line of work is to study these families via other formulations.\\

In this work we present results on the asymmetric representatives formulation \cite{Campelo08} for vertex coloring problems.
In Section \ref{sec:formulacion} we introduce this formulation in detail and we derive from it a more compact formulation by eliminating redundant variables and equalities. 
We also adapt this formulation to solve the precoloring extension and the max coloring problem.
In Section \ref{sec:genresults} we show that the vertex coloring polytope for a given graph $G$ obtained with the new compact formulation can be interpreted as the stable set polytope of another graph $\REPORD{} G$, obtained from $G$.
It turns out that $\REPORD{}{G}$ is the same graph $\tilde{G}$ presented by Cornaz and Jost in \cite{Cornaz} and analyzed afterwards in \cite{BONOMOCornaz}.
The remaining of Section \ref{sec:genresults} is devoted to study $\REPORD{}{G}$ for different classes of graphs. 
In particular, we analyze graphs whose complements have no \emph{triangle}, \emph{paw} or \emph{\kite{}} as a subgraph (see descriptions on Section \ref{sec:genresults}).
From this analysis, we derive complete characterizations for the corresponding vertex coloring polytope based on known complete characterizations for the stable set polytope.
Finally, in Section \ref{sec:finremarks} we draw some conclusions and we depict potential lines of future work.

\section{The asymmetric representatives formulation}\label{sec:formulacion}

Given a graph $G=(V,E)$, we denote by $\comp{G} = (V,\comp{E})$ the complement of $G$. 
A graph $G'=(V',E')$ is a \emph{subgraph} of $G$ if $V' \subseteq V$ and $E' \subseteq E$, and it is an \emph{induced subgraph} if $V' \subseteq V$ and $E' = \{uv \in E: u,v \in V'\}$.
We call $G'$ a \emph{spanning} subgraph if $V' = V$.
If $W \subseteq V$, then $G[W]$ is the subgraph of $G$ induced by the vertex subset $W$. 
For a vertex $u\in V$, the \emph{neighborhood} (resp. \emph{non-neighborhood}) of $u$ is $\neig(u) = \{v \in V : uv \in E\}$ (resp. $\ane{u} = \{v \in V : uv \in \comp E\}$).

\paragraph{The classical vertex coloring problem}
Given a graph $G=(V,E)$ and a total order $\prec$ on the set of vertices, the non-neighborhood of a vertex $u \in V$ can be partitioned into its \emph{lower non-neighborhood} $\nane{u} = \{v \in \ane u : v \prec u\}$ and its \emph{upper non-neighborhood} $\pane{u} = \{v \in \ane u : u \prec v\}$. 
In the asymmetric representatives formulation for the classical vertex coloring problem, for each $u \in V$ there is a binary variable $x_{uu}$ stating whether $u$ is the representative of its own color class or not. 
Additionally, for each $v \in \pane u$, a binary variable $x_{uv}$ states whether vertex $u$ is the representative of the color class assigned to vertex $v$ or not. 
Throughout this paper, we use $\X{u,W}$ and $\X{W,u}$ as shortcuts for $\sum_{v \in W} x_{uv}$ and $\sum_{v \in W} x_{vu}$, respectively.
With these definitions, a vector satisfying
		\begin{align}
			x_{uu} + \X{\nane u, u} &= 1	& \forall u \in V, \label{eq:reps_1} \\
			\X{u,K}		&\leq  x_{uu}	&\forall u \in V,~\forall K \subseteq \pane u,~K \text{ being a clique } \label{eq:reps_2}\\
			x_{uu}, x_{uv}												&\in \{0,1\}	&\forall u \in V,~\forall v \in \pane u, \label{eq:reps_bin}			
		\end{align}
represents a proper coloring of $G$ \cite{Campelo08}. Indeed, constraints (\ref{eq:reps_1}) assert that every vertex is represented by one of its lower non-neighbors or by itself. Constraints (\ref{eq:reps_2}) prevent a vertex $u$ to represent two upper neighbors if there is an edge between them. Furthermore, constraints (\ref{eq:reps_2}) allow $u$ to represent a vertex only if $u$ is a representative itself. 
The classical vertex coloring problem is solved by minimizing the number of representatives, i.e., the objective function is the sum of the variables $x_{uu}$, for $u\in V$.
Next, we analyze this formulation further to obtain a more compact one.\\

Note that in the above formulation, every variable $x_{uu}$ can be defined by equation (\ref{eq:reps_1}) to be $x_{uu} = 1 - \X{\nane u, u}$.
Hence, these variables can be eliminated from the formulation by rewriting constraints (\ref{eq:reps_2}) as
	\begin{align}
		\X{\nane u, u}  + \X{u, K} & \leq  1 &\forall u \in V,~\forall K \subseteq \pane u\text{, $K$ being a clique}  \label{eq:reps_stable}
	\end{align}
{(note that for vertices $u$ having $\pane{u} = \emptyset$ a constraint (\ref{eq:reps_stable}) with $K=\emptyset$ should be present in the formulation).
With this reformulation, the objective function must be rewritten as 
\begin{equation*}
\min \sum_{u\in V} x_{uu} = \min \sum_{u\in V} \Big(1 - \X{\nane u, u}\Big) = \max \sum_{u\in V} \X{\nane u, u} - n, 
 \end{equation*}
 which in turn is asking to maximize the number of vertices being represented, thus minimizing the number of distinct representatives.}\\

{\paragraph{The max-coloring problem} Although the formulation above corresponds to the classical vertex coloring problem, it can be adapted to solve the max-coloring problem also by using a particular ordering on the vertex set of the graph. 
As stated, given an ordering on the vertex set, the representative of a color class is the vertex with the smallest index within the class. 
Hence, by using the weight function to give the vertices a non-increasing order, the representative of a color class will be a vertex with maximum weight over the class.
The objective function for this formulation is similar to the one above but each variable $x_{wu}$ is multiplied by the weight of vertex $u$, as the following equation shows:

\begin{align*}
\min \sum_{u\in V} \omega(u)\cdot x_{uu} &= \min \sum_{u\in V} \omega(u)\cdot (1 - \X{\nane u, u})\\
			&= \max \sum_{u\in V} \omega(u)\cdot \X{\nane u, u} - \sum_{u\in V} \omega(u)
\end{align*}
where $\omega(u)$ is the weight for vertex $u$.}\\

\paragraph{The precoloring extension problem} 
Constraints (\ref{eq:reps_bin})-(\ref{eq:reps_stable}) describe proper colorings of a graph $G$. We present here a formulation using the same encoding to model the precoloring extension problem, though using a specific ordering on the vertices of the graph. 
Given a partial assignment \mbox{$\PC:V'\to\N$}, for some $V'\subseteq V$, let $\PC_c := \{ v \in V' : \PC(v) = c\}$.
We say that an ordering $\prec$ is \emph{consistent} with $\PC$ if for each vertex $w \in V\setminus V'$ and each non-empty $\PC_c$, there is at least one vertex $v \in \PC_c$ such that $v \prec w$.
For a consistent ordering $\prec$ we define 
$$
\text{rep}(v) = \min_{\prec}\{u \in \PC_{\PC(v)}\}
$$
for every $v \in V'$.
With these definitions, the formulation for the precoloring extension problem is given by (\ref{eq:reps_bin})-(\ref{eq:reps_stable}) and the following constraints
\begin{align}
x_{uv} &= 1		 &\forall uv \in \comp E \text{ such that }  u = \text{rep}(v). \label{eq:reps_preext}
\end{align}
Function $\text{rep}(v)$ identifies the representative for each class of precolored vertices for the given consistent ordering and constraints (\ref{eq:reps_preext}) assert each precolored vertex to be represented by the proper precolored representative.
On the other hand, non-precolored vertices are not constrained to any particular representative.
It is worth to note that a consistent ordering can be trivially obtained for any $\PC$, e.g., by starting the ordering with one vertex from each non-empty $\PC_c$.\\


Given a graph $G=(V,E)$, an ordering $\prec$ on $V$ and a partial assignment \mbox{$\PC:V\to\N$}, we define $\COL{G}$ (resp. $\PREEXT{G,\PC}$) to be the convex hull of the points $x\in\mathbb{R}^{|\comp E|}$ satisfying constraints (\ref{eq:reps_bin})-(\ref{eq:reps_stable}) (resp. (\ref{eq:reps_bin})-(\ref{eq:reps_preext})). 
If $\mathcal{G}$ is a family of graphs, then $\COL{\mathcal{G}}$ and $\PREEXT{\mathcal{G},\PC}$ denote the corresponding families of polytopes. 
We may omit the ordering $\prec$ in the previous definitions whenever it is clear from the context.

Given a family of polytopes $\mathcal{P}$, the associated {\em separation problem} takes a polytope $P \in \mathcal{P}$ and a vector $\hat y$ and asks to determine whether $\hat y$ belongs to $P$ or not, and if not, to find a hyperplane separating $\hat y$ from $P$. 
In turn, the {\em optimization problem} takes a polytope $P \in \mathcal{P}$ and a vector $c$ and asks for a vector $\hat x \in P$ maximizing the objective function $c^T\hat x$, unless $P = \emptyset$. 
Based on the ellipsoid method, Gr\"otschel, Lov\'asz and Schrijver \cite{GLS88} proved that the separation and the optimization problems are polynomially equivalent, i.e., if one of these problems is polynomially solvable over a family of polytopes $\mathcal P$, so is the other one. 

	%
	


\begin{theorem}\label{th:preext}
	Given a family $\mathcal G$ of graphs, if the separation problem over $\COL{G}$ can be solved in polynomial time for any $G \in \mathcal G$ and any ordering $\prec$, then for any $G \in \mathcal G$ and any partial assignment $\PC$, the separation problem over $\PREEXTORD{\prec_2}{G,\PC}$ can be solved in polynomial time for any ordering $\prec_2$ consistent with $\PC$.
\end{theorem}

\begin{proof} 
Let $\prec_2$ be an ordering consistent with $\PC$  and 
let $Q = \{x \in\R^{|\comp E|}: x_{uv} = 1$ for every $u,v \in V$ such that $u = \text{rep}(v)\}$. 
We claim that
\begin{equation} \label{eq:PREigualP}
\PREEXTORD{\prec_2}{G,\PC} = \COLORD{\prec_2}{G} \cap Q.
\end{equation}
For the reverse inclusion, take a point $\hat x \in \COLORD{\prec_2}{G} \cap Q$. Since $\hat x \in \COLORD{\prec_2}{G}$, then $\hat x$ is a convex combination of colorings $x^1,\dots,x^k$ of $G$. Since $\hat x \in Q$, then $\hat x_{uv} = 1$ for every $u,v \in V$ such that $u = \text{rep}(v)$, and this implies $\hat x^j_{uv} = 1$ for every $j=1,\dots,k$. 
Therefore, $x^j\in Q$ for $j=1,\dots,k$.
Then, all these colorings belong to $\PREEXTORD{\prec_2}{G,\PC}$ and so does $\hat x$. 
This proves that $\COLORD{\prec_2}{G} \cap Q \subseteq \PREEXTORD{\prec_2}{G,\PC}$. 

For the forward inclusion, now take a point $\hat x \in \PREEXTORD{\prec_2}{G,\PC}$. 
This point is a convex combination of colorings of $G$ where every vertex $u \in V$ such that $u = \text{rep}(u)$ is the representative of its color class and $u$ represents every other vertex $v \in V$ such that $u = \text{rep}(v)$. Hence, all these colorings belong to $Q$ and so does $\hat x$, implying $\PREEXTORD{\prec_2}{G,\PC} \subseteq \COLORD{\prec_2}{G} \cap Q$. This shows that (\ref{eq:PREigualP}) holds.

Assume the separation problem over $\subseteq \COL{G}$ can be solved in polynomial time for any ordering $\prec$. Then, a point $\hat x \notin \PREEXTORD{\prec_2}{G,\PC}$ either does not belong to $\COLORD{\prec_2}{G}$ or has $\hat x_{uv} < 1$ for some $u,v \in V$ such that $u = \text{rep}(v)$.
Hence, to separate a point from $\PREEXTORD{\prec_2}{G,\PC}$ we just need to test if $x_{uv} = 1$, for all $u,v \in V$ such that $u = \text{rep}(v)$ and, if these conditions hold, separate the point (in polynomial time) from $\COLORD{\prec_2}{G}$. 
Hence, the separation problem over $\PREEXTORD{\prec_2}{G,\PC}$ can be solved in polynomial time.
\end{proof}
\\

The mentioned equivalence between polyhedral separation and optimization \cite{GLS88} yields the following corollary.

\begin{corollary}
	Given a family $\mathcal G$ of graphs, if the optimization problem over $\COL{G}$ can be solved in polynomial time for any $G \in \mathcal G$ and any ordering $\prec$, then for any $G \in \mathcal G$ and any partial assignment $\PC$, the optimization problem over $\PREEXTORD{\prec_2}{G,\PC}$ can be solved in polynomial time for any ordering $\prec_2$ consistent with $\PC$.
\end{corollary}

Finally, the results above let us reach an important conclusion about the potential of the representatives formulation.

\begin{theorem}\label{thm:reps_potencial}
	Let $\mathcal{G}$ be a family of graphs. {If either the precoloring extension problem or the max-coloring problem} is NP-complete over $\mathcal{G}$, then the optimization/separation problem over $\COL{\mathcal G}$ with an arbitrary ordering $\prec$ is NP-complete.
\end{theorem}

\begin{proof}
Suppose the optimization problem over $\COL{G}$ can be solved in polynomial time for $G \in \mathcal G$ and for any ordering $\prec$. 
{We know that by using the weight function to order the vertex set of $G$ the max-coloring problem can be solved by optimizing over $\COL{G}$.}
Furthermore, for any partial assignment $\PC$, we can optimize over $\PREEXTORD{\prec_2}{G,\PC}$ with a consistent ordering $\prec_2$ in polynomial time, thus solving the precoloring extension problem on $G$ and $\PC$. 
Both problems being polynomially solved contradicts the hypothesis.
\end{proof}
\\

Theorem \ref{thm:reps_potencial} implies that even when the classical vertex coloring problem is polynomially solvable over $\mathcal{G}$, the polytope associated with the representatives formulation cannot be subject to an elegant characterization (for an arbitrary vertex ordering) if {either the precoloring-extension problem or the max-coloring problem} is NP-complete over $\mathcal{G}$. 
This may limit the analysis scope for the representatives formulation, as these problems are known to be NP-complete for many classes of graphs. 
Nevertheless, there are also several graph classes for which these problems can be solved in polynomial time, hence this formulation may be explored in order to find nice polyhedral descriptions for the corresponding polytopes.\\

\section{General results and polyhedral characterizations}\label{sec:genresults}

A \emph{stable set} on a graph is a set of pairwise non-adjacent vertices and the \emph{stable set polytope} $STAB(H)$ of a graph $H$ is the convex hull of the characteristic vectors of the stable sets of $H$.
As shown in the previous section, the formulation of $\COL G$ is given by (\ref{eq:reps_bin})-(\ref{eq:reps_stable}). 
On the other hand, since every non-zero coefficient in inequalities (\ref{eq:reps_stable}) is 1 and so is the right hand side, then inequalities (\ref{eq:reps_bin})-(\ref{eq:reps_stable}) give also the description of the stable set polytope on some particular graph related to $G$.\\

\begin{figure}
	\centering
	\includegraphics[height=0.3\textwidth]{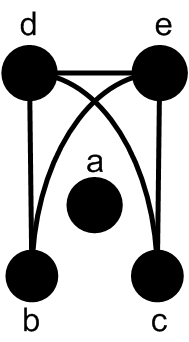} \hspace{0.10\textwidth}
	\includegraphics[height=0.3\textwidth]{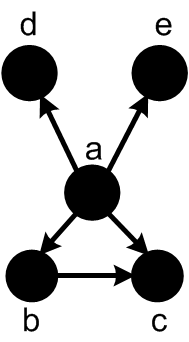} \hspace{0.05\textwidth}
	\includegraphics[height=0.3\textwidth]{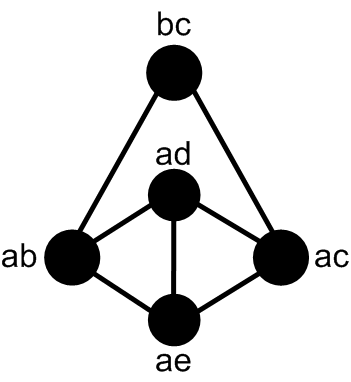}
	\caption{A graph $G$, an orientation $\prec$ of $\comp{G}$ and the resulting graph $\REP{G}$.}
	\label{fig:RG} 
\end{figure}

Define the graph $\REP G$ to have one vertex for each non-edge of $G$ (i.e., an edge of $\comp G$), and two vertices are adjacent in $\REP G$ if the corresponding non-edges in $G$ appear in the same constraint from (\ref{eq:reps_stable}). Formally the vertex set of $\REP G$ is $V(\REP G) = \comp E$ and the edge set is 
\begin{align*}
	E(\REP G) = \{(uv)(u'v') & \in \comp E\times \comp E:
						\\ & x_{uv}, x_{u'v'} \in \text{ sup($\pi, \pi_0$) } \text{for some ($\pi, \pi_0$) from (\ref{eq:reps_stable})}\},
\end{align*}
where sup$(\pi, \pi_0) = \{ x_{uv} : \pi_{uv} \neq 0\}$, i.e., the \emph{support} of the inequality given by ($\pi, \pi_0$). 
Figure \ref{fig:RG} gives an illustration of a graph $G$, an orientation $\prec$ of its complement $\comp{G}$ (an arrow from $u$ to $v$ indicates that $u \prec v$) and the resulting graph $\REP{G}$.
With these definitions, we reach the following main result:

	\begin{theorem}\label{th:coltostab}
			Given a graph $G=(V,E)$ and an ordering $\prec$ on the set of vertices, $x \in \{0,1\}^{|\comp E|}$ is the incidence vector of a proper coloring of $G$ if and only if $x$ is the incidence vector of a stable set in $\REP G$. That is, $$\COL{G} = STAB(\REP G).$$
	\end{theorem}
	
	\begin{proof}
		{
		For the forward inclusion, take a point $\hat x \in \COL{G}$. 
		Since $\hat x \in \COL{G}$, then $\hat x$ is a convex combination of colorings $\hat x^1,\dots,\hat x^k$ of $G$. 
		Call $S^j := \{x_{uv} : \hat x^j_{uv} = 1\}$, for $j=1,\dots,k$.
		Each coloring $\hat x^j$ satisfies (\ref{eq:reps_bin}) and (\ref{eq:reps_stable}), hence for each constraint $(\pi, \pi_0)$ from (\ref{eq:reps_stable}), $|S^j \cap \text{sup}(\pi, \pi_0)| \leq 1$.
		Then, by definition of $\REP{G}$, the set $S^j$ corresponds to a stable set of $\REP{G}$ and so $x^j \in STAB(\REP{G})$, for $j=1,\dots,k$. 
		Therefore, $\hat x \in STAB(\REP{G})$.
		
		For the backward inclusion, take a point $\hat x \in STAB(\REP{G})$. 
		Since $\hat x \in STAB(\REP{G})$, then $\hat x$ is a convex combination of stable sets $\hat x^1,\dots,\hat x^k$ of $\REP{G}$.
		Being a stable set, each $x^j$ satisfies (\ref{eq:reps_bin}) and (\ref{eq:reps_stable}), by definition of $\REP{G}$, and so $x^j \in \COL{G}$ for every $j=1,\dots,k$.
		Therefore, $\hat x \in \COL{G}$.}
	\end{proof}
	\\
	
	In \cite{Cornaz}, Cornaz and Jost exhibit a polynomial-time reduction from the vertex coloring problem to the maximum stable set problem. Given a graph $G$ and an acyclic orientation $D$ of $\comp G$, they construct an auxiliary graph $\tilde{G}$ such that the set of all stable sets of $\tilde{G}$ is in one-to-one correspondence with the set of all the vertex colorings from $G$. It follows from the construction of $\tilde{G}$, that this auxiliary graph coincides with $\REP{G}$.	
	{Indeed, the construction of $\REP{G}$ gives an attractive (and yet missing) interpretation of the correspondence between colorings of $G$ and stable sets of $\tilde{G}$ (i.e., $\REP{G}$), as each arc $uv$ of $\comp{G}$ expresses that $u$ represent the color assigned to $v$, thereby making the ``independence'' among the chosen arcs 	clear.
	Furthermore, Theorem \ref{th:coltostab} gives the polyhedral arguments for the mentioned correspondence between colorings and stable sets.
	}\\
	
	The \emph{line graph} $L(G)$ of a graph $G=(V,E)$ is the graph with vertex set $E$ where two vertices are linked by an edge in $L(G)$ if they correspond to two adjacent edges in $G$. 
	In \cite{Cornaz}, the graph $\tilde{G}$ is constructed by removing some edges from the line graph of $\comp G$.
	We reach the same property for $\REP G$ for our construction.\\
	
	\begin{remark} \cite{Cornaz}\label{rem:REPsubLG}
		$\REP G$ is a spanning subgraph of $L(\comp G)$.
	\end{remark}\\

	Two edges in a graph are adjacent if they share one of its endpoints. A \emph{matching} of a graph is a set of pairwise non-adjacent edges and the \emph{matching polytope} $MATCH(H)$ of a graph $H$ is the convex hull of the characteristic vectors of the matchings of $H$. The stable set polytope of the line graph of $H$ is the matching polytope of $H$. Therefore, Theorem \ref{th:coltostab} and Remark \ref{rem:REPsubLG} suggest a strong relation between $MATCH(\comp{G})$ (i.e., $STAB(L(\comp{G}))$) and $\COL{G}$.
		
	\begin{proposition}\label{prop:MatchInCol}
		For any graph $G$, we have that
		$$STAB(L(\comp G))  = MATCH(\comp G) \subseteq \COL{G} = STAB(\REP G).$$
	\end{proposition}	
	\begin{proof}
		Any matching from $\comp G$ contains at most one edge incident to a vertex $u$, for all $u \in V$. Thereby, it is trivial to check that if $x\in \{0,1\}^{|\comp E|}$ induces a matching of $\comp G$, then $x$ satisfies (\ref{eq:reps_stable}), hence it belongs to $\COL{G}$. So, every convex combination of matchings from $\comp G$ belongs to $\COL{G}$. 		
		This also follows from Remark \ref{rem:REPsubLG}, which implies that any stable set of $L(\comp G)$ is a stable set of $\REP G$, hence the inclusion between the polytopes also follows.
	\end{proof}\\

	
	The stable set polytope has been widely studied and many facet-inducing inequalities are known for this polytope. 
	Furthermore, complete characterizations are known for some graph families.	
	Due to Theorem \ref{th:coltostab}, given a family $\mathcal{G}$ of graphs, we are interested in characterizations for the family $\mathcal{\REP G} := \{ \REP G ~/~ G \in \mathcal{G} \}$, since a complete characterization of the stable set polytopes associated to $\mathcal{\REP G}$ will give us a complete characterization of the vertex coloring polytopes associated to graphs in $\mathcal{G}$.\\

\subsection{Graphs $G$ with $\alpha(G) \leq 2$}

	The \emph{stability number} $\alpha(G)$ of a graph $G$ is the cardinality of a maximum stable set of $G$.
	For graphs $G$ with $\alpha(G) \leq 2$, the vertex coloring problem can be solved in polynomial time \cite{Tuza01}. 
	
		\begin{proposition}\label{prop:ResLG}
			Given a graph $G=(V,E)$, $\alpha(G) \leq 2$ if and only if for any ordering $\prec$, we have $\REP G = L(\comp G)$.
		\end{proposition}
	
	\begin{proof}
		Assume $\alpha(G) \leq 2$ and let $u\in V$. Then, for any ordering $\prec$ on the vertices, $\pane{u}$ induces a complete subgraph (otherwise $u$ and the endpoints of any non-edge from $\pane{u}$ would induce a $K_3$ in $\comp G$). 
		Therefore, inequality (\ref{eq:reps_stable}) with the clique $\pane{u}$ implies that for every $v,w \in \nane{u} \cup \pane{u}$, vertices $uv, uw \in V(\REP G)$ are adjacent in $\REP G$.
		Thus, every pair of edges of $\comp G$ sharing $u$ as an endpoint corresponds to adjacent vertices from $\REP G$ and so, since $u$ is any arbitrary vertex, $L(\comp G)$ is a subgraph of $\REP G$. 
		Hence, by Remark \ref{rem:REPsubLG}, $\REP G = L(\comp G)$.
		
		Assume now that $\REP G = L(\comp G)$ for any ordering $\prec$ and assume that $\alpha(G) > 2$. Let $\{u,v,w\} \in V$ induce a $K_3$ in $\comp G$ and take an ordering $\prec$ such that $u \prec v \prec w$.
		With this ordering, $x_{uv}$ and $x_{uw}$ cannot simultaneously belong to the support of an inequality from (\ref{eq:reps_stable}), and so vertices $uv$ and $uw$ are not adjacent in $\REP G$, thus contradicting the fact that $\REP G = L(\comp G)$.
	\end{proof}\\

	As the matching polytope of a graph coincides with the stable set polytope of its line graph, the following is a direct consequence of Proposition \ref{prop:MatchInCol} and Proposition \ref{prop:ResLG}.
	It also gives a polyhedral argument for the well-known fact that the coloring problem over graphs $G$ with $\alpha(G) \leq 2$ can be solved by resorting to a matching problem on its complement \cite{Tuza01}.

	\begin{corollary}\label{cor:match}
		Given a graph $G$, then $\alpha(G) \leq 2$ if and only if for any ordering $\prec$, we have
		$STAB(L(\comp G))  = MATCH(\comp G) = \COL{G} = STAB(\REP G).$
	\end{corollary}

	Edmonds \cite{Edmonds65} proved that a complete characterization for the Matching Polytope of a graph $H=(V_H,E_H)$ is given by the following inequalities:
	\begin{eqnarray}			
			\sum_{v \in N(u)} x_{uv} 	& \leq 1 								&\forall u \in V_H \label{eq:match1}\\
			\sum_{uv \in E_H(S)} x_{uv} & \leq \frac{|S|-1}{2} 	&\forall S \subseteq V_H, |S| \text{ odd}\label{eq:match2}\\
			x_{uv}										& \geq 0								&\forall uv \in E_H\label{eq:match3}
		\end{eqnarray}
	
	Therefore, it follows from Corollary \ref{cor:match} that (\ref{eq:match1})-(\ref{eq:match3}) applied to the complement of a graph $G$ gives a complete characterization for $\COLORD{}{G}$ if and only if $\alpha(G) \leq 2$. 
	This result appeared in \cite{Cornaz}, however, we recast it here in our notation for the sake of completeness.
	
\begin{theorem}\label{th:cotriangle}
	\cite{Cornaz}
	Let $G=(V,E)$ be a graph with $\alpha(G) \leq 2$. For any ordering $\prec$ on the set of vertices, inequalities (\ref{eq:match1})-(\ref{eq:match3}) applied to $\comp G$ give a complete characterization for $\COL{G}$.
\end{theorem}
	
	As far as we know, the complexity of the precoloring extension problem on graphs $G$ with $\alpha(G) \leq 2$ is not known. 
	Since the separation problem for constraints (\ref{eq:match1})-(\ref{eq:match3}) can be solved in polynomial time \cite{GLS88}, the following theorem is a direct consequence of Theorem \ref{th:preext} and Theorem \ref{th:cotriangle}.
	
\begin{theorem}\label{th:preext_poly}
The precoloring extension problem can be solved in polynomial time for graphs $G$ with $\alpha(G) \leq 2$.
\end{theorem}

Another interesting conclusion arises from the reinterpretation of the \emph{odd set inequalities} (\ref{eq:match2}), in terms of the coloring polytope. 
These inequalities result to be valid for $\COLORD{}{G}$ even for a general graph $G$, whenever $\alpha(G[S]) \leq 2$. 
In fact, in this case, the resulting inequality is dominated by (or equal to) the following reformulation of the \emph{internal inequalities} from \cite{Campelo08}, as $\chi(\comp G[S]) \geq \frac{|S|+1}{2}$:

\begin{theorem}\cite{Campelo08}
		Given a graph $G=(V,E)$, an ordering $\prec$ on the set of vertices and a subset $S \subseteq V$, the \emph{internal inequality} defined as
		\begin{equation}\label{eq:dv_inner}
			\sum_{u \in S} \X{\nane{u} \cap S,u} \leq |S| - \chi(G[S])
		\end{equation}
		is valid for $\COL{G}$. In addition, (\ref{eq:dv_inner}) is facet-defining if $S$ induces an odd hole or an odd antihole.
\end{theorem}

%
%
	%


While sufficient conditions for (\ref{eq:dv_inner}) to be facet-defining are given, a complete characterization for these class of facets is missing. 
Although we fail in finding such a characterization in general, we expand below the cases for which inequalities (\ref{eq:dv_inner}) are facet-defining.\\

Given a graph $G=(V,E)$, a \emph{cut vertex} $u\in V$ is a vertex such that $G-u$ is not connected. 
$G$ is \emph{2-connected} if it has no cut vertex and it is \emph{hypomatchable} if $G - u$ has a perfect matching for every $u\in V$. 
Edmonds and Pulleyblanck \cite{EdmondsPulley} gave a characterization for the constraints (\ref{eq:match2}) defining facets of the matching polytope; namely, those for which the set $S$ induces a 2-connected hypomatchable subgraph.

\begin{theorem}\label{th:facets}
	Given a graph $G=(V,E)$ and a subset $S \subseteq V$, if $\alpha(G[S]) \leq 2$ and $\comp G[S]$ is a 2-connected hypomatchable graph, then the internal inequality (\ref{eq:dv_inner}) defines a facet of $\COL{G}$.
\end{theorem}

\begin{proof}
	Since $\comp G[S]$ is a 2-connected hypomatchable graph, then $|S|$ is odd and $\chi(G[S]) \leq \frac{|S|+1}{2}$. As $\alpha(G[S]) \leq 2$, then $\chi(G[S]) = \frac{|S|+1}{2}$. Hence, (\ref{eq:dv_inner}) is the corresponding constraint (\ref{eq:match2}).
	As (\ref{eq:match2}) defines a facet of $MATCH(\comp G)$ (which is full-dimensional) and is valid for $\COL G$, then Proposition \ref{prop:MatchInCol} implies that (\ref{eq:dv_inner}) defines a facet of $\COL G$.
%
\end{proof}

\subsection{Complements of graphs with no paw as a subgraph}

\begin{figure}
	\centering
	\includegraphics[scale=0.8]{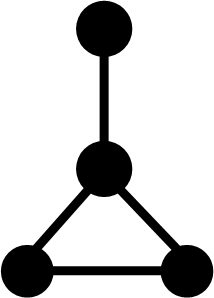}\hspace{1cm}
	\includegraphics[scale=0.8]{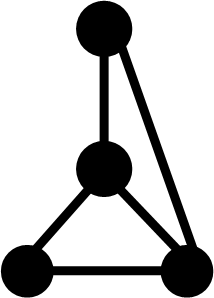} \hspace{1cm}
	\includegraphics[scale=0.8]{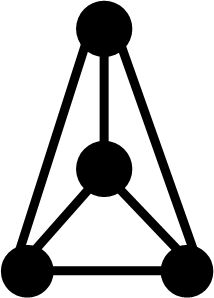} 
	\caption{The paw, the diamond and the $K_4$ graphs.}%
	\label{fig:paw}%
\end{figure}

A \emph{paw} is the graph obtained from $K_3$ by adding a fourth vertex with degree 1.
Having found a characterization for graphs $G$ with $\alpha(G) \leq 2$ (i.e., for $\comp{G}$ without $K_3$ as a subgraph), we move now to a superclass of this family. 
In this section we analyze the class of graphs such that their complements do not contain a paw as a subgraph. 
In terms of forbidden induced subgraphs, this is the class of the complements of \{$K_4$,diamond,paw\}-free graphs, were a \emph{diamond} is obtained by adding an edge to the paw. Figure \ref{fig:paw} illustrates these forbidden structures. 
We call graphs in this family to be co-\{$K_4$, diamond, paw\}-free.\\

	The complete join of two graphs $G_1=(V_1, E_1)$ and $G_2=(V_2, E_2)$ is the graph $G_1 \star G_2$ with vertex set $V_1 \cup V_2$ and edge set $E_1 \cup E_2 \cup (V_1 \times V_2)$. 
	As a general remark, it is worth to note that when a graph $G$ is the complete join of some graphs $G_1,\dots, G_k$, then the set of variables describing $\COL{G}$ is the disjoint union of the sets of variables for $\COL{G_1}, \dots, \COL{G_k}$, as there are no anti-edges among the components.
	Moreover, every constraint (\ref{eq:reps_stable}) for $G$ uses variables restricted to a unique component $G_i$, for $i \in \{1,\dots, k\}$, and so 
	$$
	\COL{G} = \{ (X^1, \dots, X^k) : X^i \in \COL{G_i}, \text{ for $i =1,\dots, k$}  \}.
	$$
	Therefore, a complete characterization for $\COL{G}$ is given by the complete characterizations for each $\COL{G_i}$.
	
	\begin{proposition} \label{prop:join}
		Let $G$ be the complete join of some graphs $G_1,\dots, G_k$ such that $\COL{G_i} = \{ x : A^ix \leq b^i \}$, for every $i = 1,\dots, k$. Then $\COL{G} = \{ (x^1,\dots, x^k) : A^ix^i \leq b^i \text{, for every $i=1,\dots, k$} \}.$
	\end{proposition}

We give next a structural characterization for the family of co-\{$K_4$, diamond, paw\}-free graphs. 

\begin{proposition}\label{prop:copaw_char}
A co-\{$K_4$, diamond, paw\}-free graph $G$ is the complete join of disjoint stable sets of size 3 and a subgraph $G'$ of $G$ with $\alpha(G') \leq 2$.
\end{proposition}

\begin{proof}
If $G=(V,E)$ has $\alpha(G) \leq 2$, then there is nothing to prove.
Thus, assume that $G$ contains a stable set $S$ of size greater than 2. As $\comp G$ has no $K_4$, then $|S| = 3$. 
For any vertex $v \in V \setminus S$, there must be an edge $uv \in E$, for every $u\in S$, otherwise $S \cup \{v\}$ would induce a paw, a diamond or a $K_4$ in $\comp G$, thus a contradiction.
Therefore, $G$ is the complete join of $S$ and $G[V\setminus S]$.
Since this applies to any stable set $S$ of size greater than 2 in $G$, the assertion follows.
\end{proof}\\

\begin{figure}
	\centering
	\begin{tabular}{r >{\centering\arraybackslash}m{2.5in}}
	$\comp{G} : $ & \includegraphics[width=.5\textwidth]{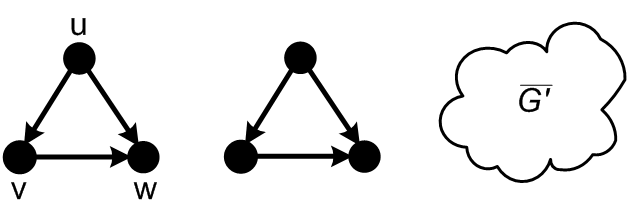}\\
	$\REP{G} : $ & \includegraphics[width=.5\textwidth]{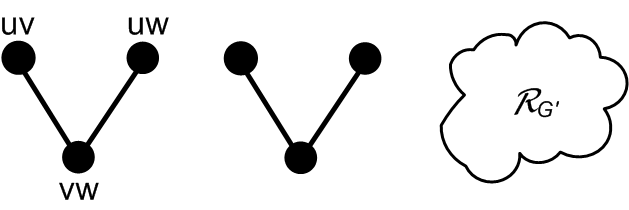}\\
	$H^{\prec}_{G} : $ & \includegraphics[width=.5\textwidth]{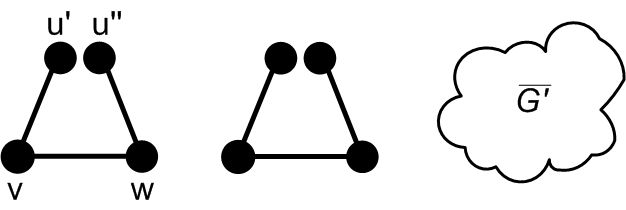}		
	\end{tabular}
	\caption{Graphs $\comp{G}$, $\REP{G}$ and $H^{\prec}_G$ for a co-\{$K_4$, diamond, paw\}-free graph $G$.}\label{fig:HG}
\end{figure}

It is easy to see that if a graph $S$ is a stable set $\{u,v,w\}$, with $u \prec v \prec w$, then $\REP{S}$ is a path of length 3, which in turn is the line graph of a path $P=\{u', v, w, u''\}$ with edges $u'v, vw$ and $wu''$. 
Hence, inequalities (\ref{eq:match1})-(\ref{eq:match3}) applied to $P$ give a complete characterization for $MATCH(P) = STAB(L(P)) = STAB(\REP{S}) = \COL{S}$.
With this fact, along with Proposition \ref{prop:join} and Proposition \ref{prop:copaw_char} and Theorem \ref{th:cotriangle}, we obtain a complete characterization for $\COL{G}$, when $G$ is a co-\{$K_4$, diamond, paw\}-free graph.
Let $G = G' \star S_1 \star \dots \star S_k$, with $\alpha(G') \leq 2$ and where each $S_i$ induces a stable set of size 3. 
Let $S_i=\{u_i,v_i,w_i\}$, where $u_i \prec v_i \prec w_i$, and define $P_i$ to be a path $\{u'_i, v_i, w_i, u''_i\}$ with edges $u'_iv_i, v_iw_i$ and $w_iu''_i$.
Define $H^{\prec}_{G}$ to be the disjoint union of $\comp{G'}$ and each of the $P_i$.
A trivial observation is that $\REP G = L(H^{\prec}_G)$ (Figure \ref{fig:HG} illustrates these constructions).

\newcommand{\GP}[1][\comp{G}]{S^{\prec}_{#1}}


\begin{corollary}\label{th:nocopaw_description}
	Let $G$ be a co-\{$K_4$, diamond, paw\}-free and $\REP{G} = L(H_G)$. 
	Then, inequalities (\ref{eq:match1})-(\ref{eq:match3}) applied to $H_G$ give a complete characterization for $\COL{G}$. 
\end{corollary}

\begin{corollary}\label{th:preext_poly_nocopaw}
The precoloring extension problem can be solved in polynomial time for co-\{$K_4$,diamond,paw\}-free graphs.
\end{corollary}


\subsection{Graphs $G$ with quasi-line $\REP{G}$}

We call a graph to be a \emph{\kite{}} if it is obtained by connecting a new vertex to the vertex of degree 1 in the paw.
A \emph{claw} is a graph formed by a vertex with three neighbors of degree one.
The graph $\tilde{G}$ introduced by Cornaz and Jost was analyzed afterwards in \cite{BONOMOCornaz} in search for structural characterizations, one of the main results being the following theorem.

\begin{theorem}[\cite{BONOMOCornaz}]\label{th:flavia}
	For any graph $G$, $\tilde{G}$ (i.e., $\REP{G}$) is a claw-free graph for any ordering $\prec$ if and only if $\comp{G}$ does not contain a \kite{} as a subgraph.
\end{theorem}

In this section we review the class of graphs $G$ such that $\comp{G}$ does not contain a \kite{} as a subgraph and we give a more precise result than Theorem \ref{th:flavia}, with some further implications.
Recall that a graph is called to be \emph{quasi-line} if the neighborhood of every vertex can be partitioned into two cliques. 
Line graphs are a subclass of quasi-line graphs which in turn are a subclass of claw-free graphs.
	
\begin{theorem}\label{th:cokite}
		For any graph $G$, $\REP G$ is a quasi-line graph for any ordering $\prec$ if and only if $\comp{G}$ does not contain a \kite{} as a subgraph.
\end{theorem}


\begin{proof}
If $\REP{G}$ is quasi-line for any ordering $\prec$, then it is claw-free and so by Theorem \ref{th:flavia}, $\comp{G}$ does not have a \kite{} as a subgraph.\\
		
Assume now that $\comp{G}$ does not have a \kite{} as subgraph.
By the construction of $\REP G$, the set of neighbors of a vertex $uv \in V(\REP G)$ is given by all the variables appearing in an inequality (\ref{eq:reps_stable}) along with $x_{uv}$. In that sense, variable $x_{uv}$ appears in 
\begin{enumerate}[(i)]
	\item every inequality (\ref{eq:reps_stable}) for vertex $v$, and 
	\item those inequalities (\ref{eq:reps_stable}) for vertex $u$ such that $v$ belongs to the corresponding clique $K \subseteq \pane{u}$.
\end{enumerate}
\newcommand{\WL}[1]{W^{\textsc l}_{#1}}
\newcommand{\WH}[1]{W^{\textsc h}_{#1}}
\newcommand{\Wuv}{W_{uv}^*}
	For each vertex $u \in V$, define 
	\begin{align*}
		\WL{u} = \{ wu \in V(\REP G): w \in \nane{u} \},\\
		\WH{u} =\{ uw \in V(\REP G): w \in \pane{u} \}.
	\end{align*}
	The set of neighbors of $uv$ arising from (i) is $\WL{v} \cup \WH{v}\setminus \{uv\}$, while the set arising from (ii) is $\WL{u} \cup \Wuv$, where $\Wuv := \{uw \in \WH{u} : vw \in E\}$. Figure \ref{fig:quasiline} gives an illustration of these sets, where thick lines (resp. thin arrows) represent edges (resp. non-edges) of $G$.	
	We will prove that $\REP G$ is a quasi-line graph by proving that for every $uv \in V(\REP G)$, the following holds:
	\begin{enumerate}[(a)]
		\item If $\WL{v} \cup \WH{v} \neq \{uv\}$, then $\REP{G}[\WL{u} \cup \Wuv]$ is a clique.
		\item If $\WL{u} \cup \Wuv \neq \emptyset$, then $\REP{G}[\WL{v} \cup \WH{v}]$ is clique.
		\item If $\WL{v} \cup \WH{v} = \{uv\}$, then $\REP{G}[\WL{u} \cup \Wuv]$ can be partitioned in two cliques.
		\item If $\WL{u} \cup \Wuv = \emptyset$, then $\REP{G}[\WL{v} \cup \WH{v}]$ can be partitioned in two cliques.
	\end{enumerate}
	
	\begin{figure}
		\centering
		\includegraphics[width=0.8\linewidth]{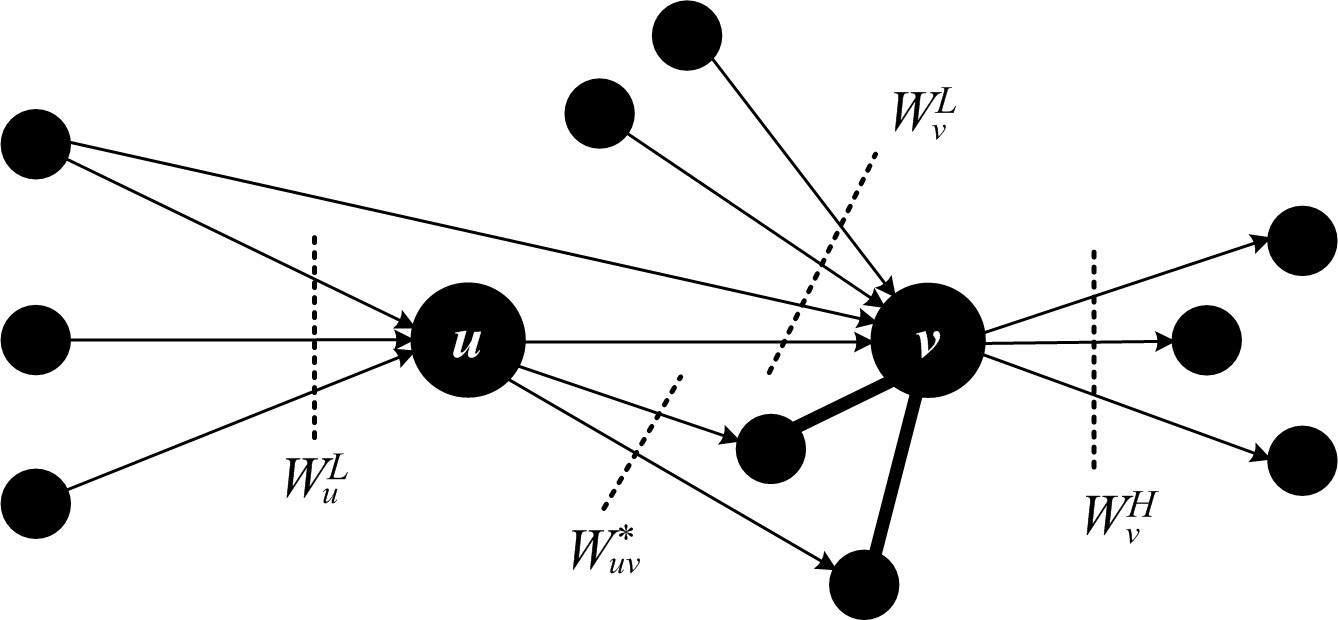}
		\caption{The neighborhood of the non-edge $uv$ in $\REP G$.}
		\label{fig:quasiline}
	\end{figure}
		
	\begin{enumerate}[(a)]
		\item 
		Let $xv \in \WL{v}$ or $vx \in \WH{v}$, with $x \neq u$. In any case, if there are $uw,uz\in \Wuv$ such that $wz\in \comp E$, then $\{x,v,u,w,z\}$ gives a \kite{} in $\comp G$. Since this cannot occur, the set $\{w \in V : vw \in \Wuv\}$ is a clique in $G$, and inequality (\ref{eq:reps_stable}) for $u$ and this clique groups every variable with subindex belonging to $\WL{u} \cup \Wuv$, and so this set induces a complete subgraph of $\REP G$.\\
		
		\item 
		Let $xu \in \WL{u}$ or $ux \cup \Wuv$. In any case, if there are $vw,vz \in \WH{v}$ such that $wz\in \comp{E}$, then $\{x,u,v,w,z\}$ gives a \kite{} in $\comp G$. As this cannot be the case, $\pane{v}$ is a clique in $G$. Hence, inequality (\ref{eq:reps_stable}) for $v$ and the clique $\pane{v}$ (or the empty set) groups every variable with subindex belonging to $\WL{v} \cup \WH{v}$, and so this set induces a complete subgraph of $\REP G$.\\
		
		\item 
		Note that $\REP{G}[\WL{u}]$ is a clique and every vertex from $\WL{u}$ is connected to every vertex from $\Wuv$, in $\REP{G}$. So, if $\REP{G}[\WL{u} \cup \Wuv]$ cannot be partitioned in two cliques then there exists an odd cycle $C$ in $\comp{\REP{G}[\Wuv]}$. Each non-edge $(uw,uz)$ of $\REP{G}[\Wuv]$ comes from a non-edge $zw\in \comp{G}$.
		If $C=\{ux,uw,uz\}$, then $\{v,u,x,w,z\}$ gives a \kite{} in $\comp G$. As this cannot be the case, assume $|C| \geq 5$ and let $ux,uw,uz$ and $uh$ be four consecutive vertices from $C$. Then $\{x,w,z,h,u\}$ gives a \kite{} in $\comp G$, hence another contradiction. Thereby, such a cycle cannot exist and so $\REP{G}[\WL{u} \cup \Wuv]$ can be partitioned in two cliques.\\
		
		\item 
		Note that $\REP{G}[\WL{v}]$ is a clique and every vertex from $\WL{v}$ is connected to every vertex from $\WH{v}$, in $\REP{G}$. As in case (c), if $\REP{G}[\WL{v} \cup \WH{v}]$ cannot be partitioned in two cliques then there exists an odd cycle $C$ in $\comp{\REP{G}[\WH{v}]}$. The rest of the proof for this case is analogous to the proof for case (c).
	\end{enumerate}

Items (a)-(c) prove that the neighborhood of any vertex $uv$ of $\REP G$ can be partitioned in two cliques. Therefore, $\REP G$ is a quasi-line graph.
\end{proof}
\\

Theorem \ref{th:cokite} gives a more precise result than Theorem \ref{th:flavia} (from \cite{BONOMOCornaz}), and also implies (along with the latter) that if $\REP{G}$ is a claw-free graph for every ordering $\prec$, then it is a quasi-line graph. With this result, we give a complete characterization for $\COL{G}$ when $\comp{G}$ has no \kite{} as a subgraph, based on the following facts.\\

Oriolo introduced the following family of valid inequalities for the stable set polytope of a graph which, in some sense, generalizes the odd-set inequalities (\ref{eq:match2}) for the matching polytope \cite{Oriolo} .
	Let $\mathcal{F} = \{K_1, \dots, K_t\}$ be a set of cliques of a graph $G=(V,E)$. Let $1 \leq p \leq t$ be an integer and $r = t$ mod $p$. Let $V_{p-1} \subseteq V$ be the set of vertices covered by exactly $(p-1)$ cliques of $\mathcal{F}$ and $V_{\geq p} \subseteq V$ the set of vertices covered by $p$ or more cliques of $\mathcal{F}$. The \emph{clique-family inequality} associated to $\mathcal{F}$ and $p$ is
	\begin{equation}\label{eq:dv_cliquefamily}
		(p-r-1) \sum_{v\in V_{p-1}} x_v + (p-r) \sum_{v\in V_{\geq p}} x_v \leq (p-r) \left\lfloor\frac{t}{p}\right\rfloor.
	\end{equation}
Oriolo proved in \cite{Oriolo} that (\ref{eq:dv_cliquefamily}) is always valid for $STAB(G)$.
The Ben Rebea's conjecture, posed in \cite{Oriolo}, says that (\ref{eq:dv_cliquefamily}) together with the non-negativity constraints and clique inequalities describe the stable set polytope of quasi-line graphs.
Some years later, this conjecture was proved to be true \cite{quasiline}.
Therefore, by Theorem \ref{th:coltostab} and Theorem \ref{th:cokite}, these inequalities give a complete characterization for $\COLORD{}{G}$, if $\comp{G}$ does not have a \kite{} as a subgraph.

\begin{theorem}
	Let $G$ be a graph such that $\comp{G}$ does not have a \kite{} as a subgraph. For any ordering $\prec$ on the set of vertices, the non-negativity inequalities together with (\ref{eq:reps_stable}) and 
	(\ref{eq:dv_cliquefamily})
	applied to $\REP G$
	give a complete characterization for $\COL{G}$.
\end{theorem}

\begin{figure}
	\centering
	\includegraphics[height=0.3\textwidth]{Counter_RG} \hspace{0.15\textwidth}
	\includegraphics[height=0.3\textwidth]{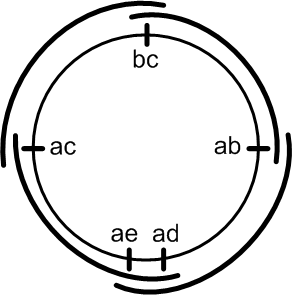} 
	\caption{The graph $\REP{G}$ from Figure \ref{fig:RG} and its FCIG representation.}
	\label{fig:FCIG} 
\end{figure}

In view of Theorem \ref{th:cokite}, one might be tempted to further strengthening the result by proving that $\REP{G}$ belongs to a certain subclass of quasi-line graphs when $\comp{G}$ has no \kite{} as a subgraph.
In fact, a quasi-line graph is either a \emph{semi-line} graph or a \emph{fuzzy circular interval graph} (FCIG), by Chudnovsky and Seymour \cite{ChudnovskySeymour2004}.
The latter graphs are defined as follows:
\begin{quote}
	\emph{Let ${\cal C}$ be a circle, ${\cal I}$ a collection of intervals in ${\cal C}$ without proper containments and common endpoints, and $V$ a multiset of points in ${\cal C}$.
	The \textit{fuzzy circular interval graph} $G(V,{\cal I})$ has node set $V$ and two nodes are adjacent if both belong to one interval $I \in {\cal I}$, where edges between different endpoints of the same interval may be omitted.}
\end{quote}
Semi-line graphs are either line graphs or quasi-line graphs without a representation as a FCIG. 
Due to a result of Chudnovsky and Seymour \cite{ChudnovskySeymour2004}, the stable set polytope of semi-line graphs requires only 0/1-valued facets.
Therefore, it is interesting to check whether all quasi-line graphs $\REP{G}$ are eventually semi-line graphs or even line graphs.
Unfortunately, there is the following counterexample. 
The graph $G$ shown in Figure \ref{fig:RG} from Section \ref{sec:genresults} is such that $\comp{G}$ has no kite as a subgraph, hence the resulting $\REP{G}$ is a  quasi-line graph.
However, $\REP{G}$ is neither a line graph (since it is one of the known forbidden induced subgraphs for line graphs \cite{Beineke}) nor a semi-line graph (since it has a representation as a FCIG, depicted here in Figure \ref{fig:FCIG}).



\section{Final remarks} \label{sec:finremarks}

The main goal of our work is to study the polytopes associated to different integer programming formulations for vertex coloring problems on particular families of graphs. 
For those cases where these problems can be polynomially solved, we pretend to find complete characterizations for the polytopes associated to at least one formulation. 
As an additional goal, we aim at extending our results in order to find such characterizations for polytopes associated to open problems, proving by this means that these variants of the vertex coloring problem can be solved in polynomial time; this is the case, for example of Theorem \ref{th:preext_poly}.
It is also expected that from this kind of studies several intermediate results arise, such as new insights in known formulations with practical implications (Theorem \ref{th:facets} being an example of this).\\

In this work we explored polytopes arising from the asymmetric representatives formulation presented in \cite{Campelo08}. Theorem~\ref{thm:reps_potencial} implies that polytopes arising from this formulation will not yield complete characterizations for those families of graphs for which the precoloring extension problem is NP-complete, unless P=NP. 
Nevertheless, there are several graph classes for which the precoloring extension problem is known to be polynomially solvable, hence we explored this formulation in order to find nice polyhedral descriptions for the corresponding polytopes.\\

The fundamental result given by Theorem \ref{th:coltostab} allows us to understand these polytopes by using available knowledge on the stable set polytope. 
This makes it possible to deduce complete characterizations for $\COL{G}$ from characterizations for $STAB$, but moreover, this shows that every family of valid inequalities known for $STAB$ gives a family of valid inequalities for $\COL{G}$.\\

As further results, we analyzed the vertex coloring polytope arising from the asymmetric representatives formulation for different classes of graphs.
In particular, we studied graphs whose complements have no \emph{triangle}, \emph{paw} or \emph{\kite{}} as a subgraph (see descriptions on Section \ref{sec:genresults}) and we derive complete characterizations for the corresponding vertex coloring polytope.\\

As a potential future line of work with the studied formulation, more families of graphs can be analyzed in search for more characterizations and new families of valid inequalities for $\COL{G}$. Also, other formulations for the vertex coloring problem are worth to be studied with the same goal, for families of graphs for which some vertex coloring problems can be solved in polynomial time.\\

\paragraph{Acknowlodgement} Part of this work was done while Diego Delle Donne was visiting LIMOS with the help of the grant \emph{``Estad\'ias cortas de doctorado en ciencia y tecnolog\'ia para profesionales argentinos en La Rep\'ublica Francesa''}, by the BEC.AR program (Argentina).

\end{document}